\documentclass[12pt]{article}
\usepackage{amsmath,amsfonts,calrsfs,amssymb,color}

\newtheorem{Theorem}{Theorem}[section]
\newtheorem{Definition}[Theorem]{Definition}
\newtheorem{Proposition}[Theorem]{Proposition}
\newtheorem{Lemma}[Theorem]{Lemma}
\newtheorem{Corollary}[Theorem]{Corollary}
\newtheorem{Remark}[Theorem]{Remark}

\newtheorem{Hypothesis}[Theorem]{Hypothesis}
\makeatletter

\@addtoreset{equation}{section}

\makeatother

\def\R{\mathbb R}
\def\N{\mathbb N}
\def\E{\mathbb E}
\def\P{\mathbb P}
 
\def\ds{\displaystyle}

\begin{document}

\title{Stochastic porous media equation and self-organized criticality}
  \author{Viorel Barbu,\\
   Institute of Mathematics ``Octav Mayer'', Iasi, Romania
,\\
Giuseppe Da Prato,\\
 Scuola Normale Superiore
di Pisa, Italy\\
 and \\
Michael R\"ockner \\
Faculty of Mathematics, University of Bielefeld, Germany\\ and\\
Department of Mathematics and Statistics, Purdue University,\\  U. S. A.}

\maketitle
{\bf Abstract}. The existence and    uniqueness of nonnegative strong solutions for 
  stochastic porous media equations with noncoercive monotone diffusivity function and Wiener forcing term is proven. The finite time extinction of solutions with high probability is also proven in $1$-$D$. The results are relevant for self-organized \textcolor{black}{critical} behaviour of stochastic nonlinear diffusion equations with critical states.



\section{Introduction}   
The phenomenon of self-organized criticality  is widely studied in Physics from different perspectives. (We refer to \cite{BJW},\cite{BJ}, \cite{Jensen}, \cite{Turco}, \cite{11}, \cite{12}, \cite{13}, \cite{14}, \cite{15} \cite{16}, \cite{17},\cite{18},\cite{19} for various studies). Roughly  speaking it is the property of systems  to have a critical point as attractor.

The purpose of this paper is to analyze this phenomenon in the frame work of stochastic evolution equations. To the best of our knowledge this is the first  time this is done in the presence of a stochastic force and in such generality in a mathematically strict way.
Let us introduce our framework.

 Let $\mathcal O$ be an open bounded domain of $\R^d, d=1,2,3,$ with smooth boundary  $\partial\mathcal O.$
We shall study here the nonlinear stochastic diffusion equation,
\begin{equation}
\label{e1.1}
\left\{\begin{array}{l}
dX(t)-\Delta\Psi(X(t))dt\ni \sigma(X(t))dW(t),\quad \quad\mbox{\rm in}\;(0,\infty)\times  \mathcal O,\\
\\
\Psi(X(t))\ni 0,\quad\mbox{\rm on}\; (0,\infty)\times \partial \mathcal O, \\
\\
X(0,x)=x\quad\mbox{\rm on}\; \mathcal O,
\end{array}\right.
\end{equation}
 where $x$ is an initial datum and $\Psi:\R\to 2^\R$ is a maximal monotone (possibly multivalued) graph with polynomial growth and random forcing term
$$
\sigma(X)dW=\sum_{k=1}^\infty \mu_kX d\beta_k\;e_k,\quad t\ge 0,
$$
where   $\{e_k\}$ is an orthonormal basis in $L^2(\mathcal O)$,  $\{\mu_k\}$ is a sequence of positive numbers and    $\{\beta_k\}$ 
a sequence of  independent   standard Brownian
motions on a filtered probability space
$(\Omega,\mathcal F,\{\mathcal F_t\}_{t\ge 0},\P)$.  

We note that $\sigma(X)$ is defined by
$$
\sigma(X)h=\sum_{k=1}^\infty \mu_kX \langle h,e_k   \rangle_2 e_k,\quad\forall\;h\in L^2(\mathcal O),
$$
where $ \langle \cdot,\cdot   \rangle_2$ is the scalar product in $L^2(\mathcal O)$.

 The equation models the dynamics of flows in porous media and more generally the phase transition (including melting and solidification processes) in the presence of a random forcing term $\sigma(X)dW$.\bigskip
 
 Existence for stochastic equations of the form \eqref{e1.1} with additive and multiplicative noise   was studied in \cite{3} under the main assumption
 that  $\Psi$ is monotonically increasing, continuous and such that
\begin{equation}
\label{e1.2}
\left\{\begin{array}{l}
\Psi(0)=0,\;\Psi'(r)\le \alpha_1 |r|^{m-1}+\alpha_2,\quad\forall\;r\in \R,\\
\\
\ds\int_0^r\Psi(s)ds\ge \alpha_3 |r|^{m+1}+\alpha_4,\quad\forall\;r\in \R,
\end{array}\right.
\end{equation}
where $\alpha_1\ge 0,\alpha_3>0$, $\alpha_2, \alpha_4\ge 0$ and $m\ge 1$.  (See also \cite{4} and \cite{5} for general growth conditions on  $\Psi$.)

Here we shall study equation \eqref{e1.1} under the following assumptions.

\begin{Hypothesis} 
\label{h1.1}
\begin{enumerate}
\item[(i)] $\Psi$ is a maximal monotone multivalued function from $\R$ into $\R$ such that $0\in \Psi(0)$.

\item[(ii)] There exist $C>0$ and $m\ge 1$ such that
$$
\sup\{|\theta|:\;\theta\in \Psi(r)\}\le C(1+|r|^m),\quad\forall\;r\in \R.
$$

\item[(iii)]  The sequence   $\{\mu_k\}$  is such that
$$
\sum_{k=1}^\infty\mu^2_k\lambda^2_k<+\infty,
$$
where $\lambda_k$ are the eigenvalues of    the Laplace operator $-\Delta$ in $\mathcal O$ with Dirichlet boundary conditions.

{\rm We recall  that the domain of $\Delta$ is $H^2(\mathcal O)\cap H^1_0(\mathcal O)$.}

\end{enumerate}
\end{Hypothesis} 
Since for $x\in H^{-1}(\mathcal O)$
\begin{equation}
\label{e1.3}
|xe_k|_{-1}^2\le C_1|e_k|^2_{H^2(\mathcal O)}\;|x|^2_{-1}
\le C_1\lambda_k^2|x|^2_{-1}
\end{equation}
and hence
\begin{equation}
\label{e1.4}
\|\sigma(x)\|^2_{L_2(L^2(\mathcal O),H^{-1}(\mathcal O))}=\sum_{k=1}^\infty\mu_k^2|xe_k|_{-1}^2
\le C_1\sum_{k=1}^\infty\mu_k^2\lambda_k^2|x|^2_{-1},
\end{equation}
 it follows by {\it (iii)} that $\sigma(x)\in L_2(L^2(\mathcal O),H^{-1}(\mathcal O))$ (the space of all  Hilbert-Schmidt operators from $L^2(\mathcal O)$ into $H^{-1}(\mathcal O)$) and that it is Lipschitz continuous from  $H^{-1}(\mathcal O)$ into $L_2(L^2(\mathcal O),H^{-1}(\mathcal O))$.
Under these assumptions we  shall prove that if $x\in L^p(\mathcal O), p\ge \max\{2m,4\},$ then there is a unique strong solution to equation \eqref{e1.1} which  is nonnegative if so is the initial data $x$.
With respect to the situation considered in  \cite{4}, in the present case one does not assume that the range of $\Psi$ is all of $\R$, which is quite unusual for porous media equations. Also Hypothesis \ref{h1.1}{\it (i)} allows monotonically increasing functions $\Psi$ with a finite number of discontinuities (jumps), $r_1,...,r_N$. One must, of course, fill the jumps by taking $\Psi(r_j)=[\Psi(r_j+0),\Psi(r_j-0)],\;j=r_1,...,r_N.$

It should be mentioned that several physical problems with free boundary and with phase transition can be put into this functional setting. For instance if
\begin{equation}
\label{e1.5}
\Psi(x)=\left\{\begin{array}{l}
\alpha_1(x-a),\quad\mbox{\rm for}\;x<a\\
\protect [0,\rho],\quad\mbox{\rm for}\;x=a\\
\alpha_2(x-a)+\rho,\quad\mbox{\rm for}\;x>a,
\end{array} \right.
\end{equation}
with $a,\rho,\alpha_1,\alpha_2\in (0,+\infty),$
then \eqref{e1.1}  models the phase transition in porous media or in heat conduction (Stefan problem). If $\Psi(x)=\rho$ sign $x$  where $\rho>0$ and
\begin{equation}
\label{e1.6}
\quad\mbox{\rm sign}\;x=\left\{\begin{array}{l}
\ds\frac{x}{|x|},\quad\mbox{\rm if}\;x\neq 0\\\\
\protect[-1,1],\quad\mbox{\rm if}\;x=0,
\end{array} \right.
\end{equation}
then \eqref{e1.1}  reduces to the nonlinear singular diffusion equation
$$
dX(t)-\rho\;\mbox{\rm div}\;(\delta(X(t))\nabla X(t))dt=\sigma(X(t))dW(t),
$$
where $\delta$ is the Dirac measure concentrated at the origin.

Other examples such as  the  Heavside step function
$$
H(x)=\left\{\begin{array}{l}
0,\quad\mbox{\rm if}\;x<0\\
\protect[0,1],\quad\mbox{\rm if}\;x=0\\
1,\quad\mbox{\rm if}\;x>0,
\end{array} \right.
$$
or $\Psi(x)=|x|^\alpha\;{\rm sign}\;x$ with $0<\alpha\le 1$   also satisfy Hypothesis \ref{h1.1}. 

In particular the equation
\begin{equation}
\label{e1.7}
dX(t)-\Delta (H+\lambda)(X(t)-x_c)dt=\sigma(X(t)-x_c)dW(t),
\end{equation}
where $\lambda>0 $, represents the continuous, stochastic version of the Bak, Jang, Wiesenfeld sand pile model \cite{BJW}. (See \cite{BJW} for a deterministic presentation of the model.)
This is a diffusion problem with free boundary driven by a random forcing term proportional to $X(t)-x_c$ where $x_c$ is the critical density and $X(t)$ is the density at the moment $t$.

Taking into account the numerical simulation in $1$-$D$ (see \cite{BJ}), one might expect that the time evolution of the system displays  self-organized criticality, i.e. the supercritical region $\{X(t)>x_c\}$
is  absorbed asymptotically in time by the critical one $\{X(t)=x_c\}$. Here we shall prove that, e.g. in all examples \eqref{e1.5}-\eqref{e1.7} above, this   indeed takes places with high probability under appropriate assumptions on the parameters
 and more precisely  that the supercritical region ``vanishes'' into the critical one in finite time with high probability, at least if $\mu_k=0$
for all $k\ge N+1$ for some $N\in \N$.  We emphasize that this  \textcolor{black}{is in particular true}
 when the noise is zero. \textcolor{black}{In this case one gets an explicit bound for the time when this happens} (cf. Remark \ref{r4.4} below).\bigskip

The plan of this paper is the following. The main results are  presented in Section 2
and are proven in Section 3.  In Section 4 we prove a finite time extinction type result for solutions to \eqref{e1.1} which displays a self-organized criticality behaviour.\bigskip

The following notations will be used.
$L^p(\mathcal O),\;p\ge 1,$ is the usual space of $p$-integrable functions with norm denoted by $|\cdot|_p$. The scalar product in $L^2(\mathcal O)$ and the duality induced by the pivot space $L^2(\mathcal O)$
will be denoted by $\langle \cdot, \cdot   \rangle_2$.
$H^k(\mathcal O)\subset L^2(\mathcal O),\;k=1,2,$ are the standard Sobolev spaces on $\mathcal O$, while $H^1_0(\mathcal O)$ is the subspace of $H^1(\mathcal O)$ with
zero trace on the boundary. For $p,q\in [1,+\infty]$ by  $L^q_W((0,T);L^p(\Omega;H))$ ($H$ a Hilbert space) we shall denote the space of all $q$-integrable processes $u:[0,T]\to L^p(\Omega;H)$ which are adapted to the  filtration $\{\mathcal F_t\}_{t\ge 0}.$  

By $C_W([0,T];L^2(\Omega;H))$  we shall denote the space of all $H$-valued adapted processes which are mean square continuous.
$L(H)$ denotes the space of bounded linear operators equipped with the usual norm.

In the following by $H$ we shall denote  the distribution space
$$
H=H^{-1}(\mathcal O)=(H^{1}_0(\mathcal O))'
$$
endowed with the scalar product and norm defined by
$$
\langle u,v   \rangle=\int_\mathcal OA^{-1}u(\xi)v(\xi)d\xi,\quad |u|_{-1}=\langle u,u   \rangle^{1/2},
$$
where $A=-\Delta$ with $D(A)=H^2(\mathcal O)\cap H^1_0(\mathcal O)$.

In terms of $A$ equation \eqref{e1.1} can be formally rewritten as
\begin{equation}
\label{e1.8}
\left\{\begin{array}{l}
dX(t)+A\Psi(X(t))dt\ni \sigma(X(t))dW(t),\ \\
\\
X(0,x)=x.
\end{array}\right.
\end{equation}
Its exact meaning will be precised later (see Definition \ref{d2.1} below).

It should be recalled, however, that the operator $x\to A\Psi(x)$ with the domain
\textcolor{black}{$$\{ x\in L^1( \mathcal O)\cap H^{-1}( \mathcal O): \mbox{\rm there is}\; \eta\in H^1_0( \mathcal O), \eta\in \Psi(x) \;\mbox{\rm a.e. in}\;  \mathcal O \}$$}
is maximal monotone in $H:=H^{-1}(\mathcal O)$ (see e.g. \cite{0}) and so the distribution space $H$  offers the natural functional setting for the porous media equation \eqref{e1.1}
or  its abstract form  \eqref{e1.8}. However, the general existence theory of infinite dimensional stochastic equations in Hilbert space with nonlinear maximal monotone operators (see \cite{DPZ2}, \cite{PR}) is not applicable in the present case and so a direct approach must be used.

Fnally, in this paper we use the same letter $C$ for  several different positive constants arising  in chains of estimates.

\section{Existence, uniqueness and positivity}
 
\begin{Definition}
\label{d2.1}
Let $x\in H$. An $H$-valued   continuous $\mathcal F_t$-adapted process $X=X(t,x)$ is called a solution to \eqref{e1.1} $($equivalently \eqref{e1.8}$)$ on $[0,T]$ if 
$$
X\in L^{p} (\Omega\times(0,T)\times \mathcal O)\cap L^2(0,T;L^2(\Omega,H)),\quad p\ge m,
$$
 and there exists $\eta\in 
L^{p/m} (\Omega\times(0,T)\times \mathcal O)$ such that $\P$-a.s.
\begin{equation}
\label{e2.1}
\begin{array}{lll}
\langle X(t,x),e_j\rangle_2 &=&\ds\langle x,e_j\rangle_2 +\int_0^t\int_\mathcal O\eta(s,\xi)\Delta e_j(\xi) d\xi ds\\
\\
&&\ds +\sum_{k=1}^\infty\mu_k\int_0^t\langle X(s,x)e_k,e_j\rangle_2 d\beta_k(s),\quad \forall\;j\in \N,\;t\in [0,T],
\end{array}
\end{equation}
\begin{equation}
\label{e2.2}
\eta\in \Psi(X)\;\quad\mbox{\rm a.e. in }\;\Omega\times (0,T)\times \mathcal O.
\end{equation}
\end{Definition}
Below for simplicity we often write $X(t)$ instead of $X(t,x)$.

From the stochastic point of view the solution $X$ given by Definition
\ref{d2.1} is a strong one, but from the PDE point of view it is a solution in the sense of distributions since the boundary condition  $\Psi(X)\textcolor{black}{\notin}0$ on $\partial \mathcal O$ is satisfied in a weak sense only.  

Theorem \ref{t2.2} below is the main existence result.
\begin{Theorem}
\label{t2.2}
Assume that  $d=1,2,3$ and that Hypothesis $\ref{h1.1}$ holds.  Then for each $x\in L^p( \mathcal O)$, $p\ge \max\{2m,4\}$  there is a unique solution $X\in L^\infty_W(0,T;L^p( \Omega;\mathcal O))$
to \eqref{e1.1}. Moreover, if $x$ is nonnegative a.e. in $\mathcal O$ then $\P$-a.s.
$$
X(t,x)(\xi)\ge  0, \quad\mbox{\it for a.e.}\;(t,\xi)\in (0,\infty)\times \mathcal O.
$$
\end{Theorem}

As mentioned earlier, Theorem \ref{t2.2} was proven in \cite{3} for a differentiable $\Psi$ satisfying conditions \eqref{e1.2} and for $p\ge \max\{m+1,4\}$. It should be said, however, that in contrast with what happens for coercive functions $\Psi$
arising in \cite{3}, here it seems no longer possible to extend the existence result to all $x\in H^{-1}(\mathcal O),\;x\ge 0$.

\section{Proof of Theorem \ref{t2.2}}

We shall consider the approximating equation
\begin{equation}
\label{e3.1}
\left\{\begin{array}{l}
dX_\lambda(t)+A(\Psi_\lambda(X_\lambda(t))+
\lambda X_\lambda(t))
dt=\sigma(X_\lambda(t))dW(t),\ \\
\\
X_\lambda(0,x)=x,
\end{array}\right.
\end{equation}
where $\lambda>0$ and
$$
\Psi_\lambda(x)=\frac1\lambda\;(x-(1+\lambda\Psi)^{-1}(x))\in \Psi((1+\lambda\Psi)^{-1}(x))
$$
is the  Yosida approximation of $\Psi$. We recall that $\Psi_\lambda$ is Lipschitzian and monotonically increasing and so
$x\to \Psi_\lambda(x)+ \lambda x$ is strictly monotonically increasing and bounded by $C_1(1+|x|^m)$ and $(\Psi_\lambda(x)+\lambda x)x\ge \lambda|x|^2$
for all $x\in \R$. By \cite[Theorem 2.2]{3} (applied with $m=1$),  for each $x\in H^{-1}(\mathcal O)$ equation \eqref{e3.1} has a unique solution
$$
X_\lambda\in L^{2}(\Omega\times (0,T)\times\mathcal O)\cap L^2_W(\Omega,C([0,T];H))
$$
in the sense of Definition \ref{d2.1}. Here as usual $C([0,T];H)$ is equipped with the supremum norm.
Moreover, ( see e.g. \cite[Theorem 4.2.5]{PR}) the  following It\^o formula holds
\begin{equation}
\label{e3.2}
\begin{array}{lll}
\E|X_\lambda(t)|^2_{-1}&+&\ds 2\E\int_0^t\int_\mathcal O(\Psi_\lambda(X_\lambda(s))+
\lambda X_\lambda(s))X_\lambda(s) d\xi\;ds\\
\\
&=&\ds |x|^2_{-1}+\sum_{k=1}^\infty\mu_k^2\;\E\int_0^t|X_\lambda(s)e_k|^2_{-1}ds.
\end{array} 
\end{equation}
We note that since
$$
|X_\lambda e_k|_{-1}\le C|e_k|_{H^2(\mathcal O)}|X_\lambda |_{-1}\le 
C \lambda _k|X_\lambda|_{-1},
$$
(cf. \eqref{e1.3}) we have by Hypothesis \ref{h1.1}(iii) (cf. \eqref{e1.4})
\begin{equation}
\label{e3.3}
\sum_{k=1}^\infty \mu_k^2 \E\int_0^t|X_\lambda(s) e_k|^2_{-1}ds\le C\E\int_0^t|X_\lambda(s)|^2_{-1}ds.
\end{equation}

\begin{Lemma}
\label{l3.1} 
There exists a constant $C>0$ such that for all $p\ge 2$ and all $x\in L^p(\mathcal O)$,
\begin{equation}
\label{e3.4}
\mbox{\rm ess.sup}_{t\in[0,T]}\;\E|X_\lambda(t,x)|^p_p\le \exp  \left(C \frac{p-1}{2} \right)\;
|x|^p_p,\quad \forall\;\lambda>0.
\end{equation}
\end{Lemma}
{\bf Proof}. We know from \cite[Lemma 3.4]{3} (with $m=1$) that as $\varepsilon\to 0$
\begin{equation}
\label{e3.5}
\left\{\begin{array}{l}
X_\lambda^\varepsilon\to X_\lambda\quad\mbox{\it strongly in}\;L^\infty_W(0,T;L^2(\Omega;H)),\\
\\
X_\lambda^\varepsilon\to X_\lambda\quad\mbox{\it weakly star in}\;L^\infty_W(0,T;L^p(\Omega;L^p(\mathcal O))),
\end{array}\right.
\end{equation}
where $X_\lambda^\varepsilon$ is the solution to  the approximating equation
\begin{equation}
\label{e3.6}
\left\{\begin{array}{l}
dX_\lambda^\varepsilon(t)+(A_\lambda)_ \varepsilon X_\lambda^\varepsilon(t)dt=\sigma(X_\lambda^\varepsilon(t))dW(t),\quad t\ge 0,\\
\\
X_\lambda^\varepsilon(0)=x,
\end{array}\right.
\end{equation}
where
$$
\left\{\begin{array}{l}
A_\lambda x=A(\Psi_\lambda(x)+\lambda x)=-\Delta (\Psi_\lambda(x)+\lambda x),\\
\\
D(A_\lambda)=\{x\in H\cap L^1(\mathcal O):\;\Psi_\lambda(x)+\lambda x\in H^1_0(\mathcal O)\},
\end{array}\right.
$$
and $(A_\lambda)_ \varepsilon$ is the Yosida approximation of $A_ \lambda$,
$$
(A_\lambda)_ \varepsilon=\frac1\varepsilon \;(I-(I+\varepsilon A_\lambda)^{-1}), \quad \varepsilon>0.
$$
Furthermore, by \cite[Lemma 3.2]{3} we have that
$
X_\lambda^\varepsilon\in L^2(\Omega;C([0,T];L^2(\mathcal O)).
$
As  a matter of fact the results of \cite{3} were proven for smooth nonlinear functions while $\Psi_\lambda$ is only Lipschitz; but the extension to  lipschitzian functions $\Psi$ satisfying \eqref{e1.2} is immediate. In fact, one might take a  smoother approximation of  $\Psi$,
for instance the mollifier $\Psi_\lambda*\rho_\lambda$ ($\rho_\lambda(r)=\frac1\lambda\;
\rho(\lambda/r),\rho\in C^\infty_0(\R), \rho\ge 0, \int\rho dr=1$) which still remains monotonically increasing and has all properties of $\Psi_\lambda$.

Next we apply It\^o's formula    \eqref{e3.6} for the function
$
\varphi(x)=\frac1p\;|x|^p_{p}.
$
More precisely, we first   apply It\^o's formula to $\varphi_\gamma(x)=\frac1p\;|(1+\gamma A)^{-1}x|^p_{p}$, $\gamma>0$, and then we let  $\gamma\to 0$.
We have (for details see the proof in \cite[Lemma 3.5]{3}),
\begin{equation}
\label{e3.7}
\begin{array}{l}
\ds\E\varphi(X^\varepsilon_\lambda(t))+\E\int_0^t\langle (A_\lambda)_ \varepsilon X^\varepsilon_\lambda(s),|X^\varepsilon_\lambda(s)|^{p-2}X^\varepsilon_\lambda(s)  \rangle_2ds\\
\\
=\ds\varphi(x)+\frac{p-1}2\;\sum_{k=1}^\infty\mu_k^2\E\int_0^t\int_\mathcal O|X^\varepsilon_\lambda(s)|^{p-2}|X^\varepsilon_\lambda(s)e_k|^{2}d\xi\;ds\;d\xi\\
\\
\le \ds \varphi(x)+\frac{p-1}2\;C \E\int_0^t\int_\mathcal O|X^\varepsilon_\lambda(s)|^{p}d\xi\;ds,
\end{array} 
\end{equation}
since by Sobolev embedding $|e_k|_\infty\le C \lambda_k$ for all $k\in \N$.
If $Y^\varepsilon_\lambda$ is  the solution to the equation
$$
Y^\varepsilon_\lambda-\varepsilon\Delta(\Psi_\lambda(Y^\varepsilon_\lambda)+\lambda Y^\varepsilon_\lambda)
=X^\varepsilon_\lambda, \quad \Psi_\lambda(Y^\varepsilon_\lambda)
+\lambda Y^\varepsilon_\lambda\in H^1_0(\mathcal O),
$$
then (see \cite[(3.25)]{3}) $|Y^\varepsilon_\lambda|_p\le |X^\varepsilon_\lambda|_p$
and therefore
$$
\langle (A_\lambda)_ \varepsilon X^\varepsilon_\lambda,|X^\varepsilon_\lambda|^{p-2} X^\varepsilon_\lambda \rangle_2=\frac1\varepsilon\;
\langle X^\varepsilon_\lambda-Y^\varepsilon_\lambda,|X^\varepsilon_\lambda|^{p-2} X^\varepsilon_\lambda \rangle_2\ge 0.
$$
Then by \eqref{e3.7} it follows, via Gronwall's lemma, that
$$
\E|X^\varepsilon_\lambda(t)|^p_p\le |x|^p_p\;\exp\left(C\frac{p-1}2\;  \right),
$$
where $C$ is independent of $x,\lambda$ and $t$.
Now one obtains \eqref{e3.4} by letting $\varepsilon$ tend to $0$ and taking  into account 
\eqref{e3.5}. $\Box$\bigskip

From now on let us assume that $p\ge\max\{4,2m\}$ and
$x\in L^p(\mathcal  O).$
 From Lemma \ref{l3.1} it follows that for a subsequence $\{\lambda\}\to 0$ we have
\begin{equation}
\label{e3.8}
\left\{\begin{array}{l}
X_\lambda\to X\quad\mbox{\rm weakly in}\;L^p(\Omega\times(0,T)\times\mathcal  O),\\
\hspace{20mm}\mbox{\rm and weakly star in  in}\;L^\infty(0,T;L^p(\Omega;L^p(\mathcal  O))),\\
\\
\Psi_\lambda(X_\lambda)\to \eta\quad\mbox{\rm weakly in}\;L^{p/m}(\Omega\times(0,T)\times\mathcal  O),\;\\
\\
\hspace{30mm}\mbox{\rm in particular in}\;L^{2}(\Omega\times(0,T)\times\mathcal  O),
\end{array}\right. 
\end{equation}
because by Hypothesis({\it ii}),
$$
|\Psi_\lambda(x)|\le |\Psi^0(x)|\le C(1+|x|^m),\quad \forall\;x\in \R.
$$
($\Psi^0$ is the minimal section of $\Psi$). By\eqref{e3.4} we  have   for $\lambda\to 0$
\begin{equation}
\label{e3.9}
 \lambda X_\lambda\to 0\quad\mbox{\rm strongly  in}\;L^{p}(\Omega\times(0,T)\times\mathcal  O).
\end{equation}
Clearly  $X$ and $\eta$ are adapted processes.
On the other  hand, we  have
$$
\begin{array}{l}
d(X_\lambda(t)-X_\mu(t))-\Delta(\Psi_\lambda(X_\lambda(t))-\Psi_\mu(X_\mu(t)) +\lambda X_\lambda(t)-\mu X_\mu(t))dt
\\
\\
=(\sigma(X_\lambda(t))-\sigma(X_\mu(t) ))dW(t)
\end{array}
$$
and therefore once again applying  It\^o's formula  (cf. \eqref{e3.2}) we obtain for $\alpha>0, t\in [0,T]$,
\begin{equation}
\label{e3.10}
\begin{array}{l}
\ds\frac12\;|X_\lambda(t)-X_\mu(t))|_{-1}^2 e^{-\alpha t}\\
\\
\ds+ \int_0^t\int_\mathcal O
\Big[(\Psi_\lambda(X_\lambda(s))-\Psi_\mu(X_\mu(s))\;(\lambda
\Psi_\lambda(X_\lambda(s))-\mu\Psi_\mu(X_\mu(s)))\\
\\
\ds\hspace{10mm}+(\lambda X_\lambda(s)-\mu X_\mu(s))( X_\lambda(s)- X_\mu(s))\Big]e^{-\alpha s} d\xi\;d s\\
\\
\ds\le \left(C\sum_{k=1}^\infty\mu_k^2\lambda_k^2-\frac12\;\alpha\right)\int_0^t|X_\lambda(s-X_\mu(s))|_{-1}^{2}e^{-\alpha s}\;ds +M_{\lambda,\mu}(t),\quad \forall\;\lambda,\mu>0,
\end{array} 
\end{equation}
where
$$
M_{\lambda,\mu}(t):= \int^t_0 e^{-\alpha s}\langle X_\lambda (s) - X_\mu (s),\sigma(X_\lambda (s) - X_\mu (s))dW(s)\rangle_2 
$$
 is a real local valued martingale. To derive \eqref{e3.10} we used that $x=\lambda\Psi_\lambda(x)+(1+\lambda\Psi)^{-1}(x)$ and thus for all $x,y\in \R$
$$
\begin{array}{lll}
(\Psi_\lambda(x)-\Psi_\mu(y))(x-y)&=&[\Psi_\lambda(x)-\Psi_\mu(y)][
(1+\lambda\Psi)^{-1}(x)-(1+\mu\Psi)^{-1}(y)]
\\
\\
&&+[\Psi_\lambda(x)-\Psi_\mu(y)][
\lambda\Psi_\lambda(x)-\mu\Psi_\mu(y)],
\end{array}
$$
and that the first summand on the right hand side is nonnegative
because $\Psi$ is  monotonically increasing and $\Psi_\lambda(x)\in \Psi((1+\lambda\Psi)^{-1}(x))$.
Hence for $\alpha>0$ large enough we obtain for all $\lambda,\mu\in (0,1)$ and $t\in[0,T]$
\begin{equation}
\label{e3.11}
\begin{array}{l}
\ds\frac12\;|X_\lambda(t)-X_\mu(t))|_{-1}^2 e ^{-\alpha t}\\
\\
\ds\le C \max\{\lambda,\mu\}\int_0^t\int_\mathcal O
\Big(|\Psi_\lambda(X_\lambda(s))|^2+|X_\lambda(s)|^2+|\Psi_\mu(X_\mu(s))|^2\\
\\
\hspace{10mm}+|X_\mu(s)|^2\Big)e^{-\alpha s}d\xi\;ds +M_{\lambda,\mu}(t).
\end{array} 
\end{equation}
Hence by the Burkholder-Davis-Gundy inequality (for  $p=1$) we get for all $\lambda,\mu \in (0,1),$ $r\in [0,T]$,
\begin{equation}
\label{e3.12}
\begin{array}{l}
\ds\frac12\;\E\sup_{t\in [0,r]}|X_\lambda(t)-X_\mu(t))|_{-1}^2 e ^{-\alpha t}\\
\\
\ds\le C \max\{\lambda,\mu\}\E\int_0^r\int_\mathcal O
\Big(|\Psi_\lambda(X_\lambda(s))|^2+|X_\lambda(s)|^2+|\Psi_\mu(X_\mu(s))|^2\\
\\
\ds+|X_\mu(s)|^2\Big)e^{-\alpha s}d\xi\;ds +C\E\left(\int_0^r|X_\lambda(s)-X_\mu(s)|^4_{-1}e ^{-2\alpha s}ds\right)^{1/2}.
\end{array} 
\end{equation}
But
\begin{equation}
\label{e3.13}
\begin{array}{l}
\ds \E\left(\int_0^r|X_\lambda(s)-X_\mu(s)|^4_{-1}e ^{-2\alpha s}ds\right)^{1/2}\\
\\
\ds\le \E\sup_{s\in [0,r]}|X_\lambda(s)-X_\mu(s))|_{-1} e ^{-\frac\alpha2 s } 
  \left(\int_0^r|X_\lambda(s)-X_\mu(s)|^2_{-1}e ^{-\alpha s}ds\right)^{1/2}\\
  \\
  \ds\le \frac14\;\E\sup_{s\in [0,r]}|X_\lambda(s)-X_\mu(s))|^2_{-1} e ^{-\alpha s }+C\E \int_0^r|X_\lambda(s)-X_\mu(s)|^2_{-1}e ^{-\alpha s}ds.
\end{array} 
\end{equation}
Taking into  account that by Hypothesis \ref{h1.1}(ii)
$$
|\Psi_\lambda(X_\lambda)|\le C(1+|X_\lambda|^m),\quad \forall\;\lambda>0,
$$
and that by   \eqref{e3.4} $\{X_\lambda\}$ is bounded in $L^{p}(\Omega\times(0,T)\times\mathcal  O)$
for  $p\ge \max\{4,2m\}$,  we infer by \eqref{e3.12}, \eqref{e3.13} 
and Gronwall's lemma   that
 $\{X_\lambda\}$ is a Cauchy net in $L^2(\Omega;C([0,T];H))$ Hence for $\lambda\to 0$
 \begin{equation}
\label{e3.14}
X_\lambda\to X\quad\mbox{\rm  in}\;L^2(\Omega;C([0,T];H)).
\end{equation}
In order to  complete the proof of the existence part   of Theorem \ref{t2.2} it suffices to show that
\begin{equation}
\label{e3.15}
\eta(\omega,t,\xi)\in \Psi(X(\omega,t,\xi))\quad\mbox{\rm a.e in }\;\Omega\times(0,T)\times\mathcal  O.
\end{equation}
Since the operator
$$
L^p(\Omega\times(0,T)\times\mathcal  O)\to L^{\frac{p}{m}}(\Omega\times(0,T)\times\mathcal  O)\subset L^{\frac{p}{p-1}}(\Omega\times(0,T)\times\mathcal  O),\;\; X\to \Psi(X),
$$
in the duality pair
$$
\left(L^p(\Omega\times(0,T)\times\mathcal  O),L^p(\Omega\times(0,T)\times\mathcal  O)'=L^{\frac{p}{p-1}}(\Omega\times(0,T)\times\mathcal  O)\right),
$$
is  maximal  monotone,  it suffices to show that (see e.g. \cite{0})
\begin{equation}
\label{e3.16}
\liminf_{\lambda\to 0}\E\int_0^T\int_\mathcal  O\Psi_\lambda(X_\lambda)X_\lambda d\xi  dt \le
\E\int_0^T\int_\mathcal  O\eta X d\xi  dt.
\end{equation}
To prove \eqref{e3.16} we first note  that by   \eqref{e3.2}  we have
\begin{equation}
\label{e3.17}
\begin{array}{l}
\ds \liminf_{\lambda\to 0}\E\int_0^T\int_\mathcal  O\Psi_\lambda(X_\lambda)X_\lambda 
d\xi  dt+\frac12\;\E|X(t)|^2_{-1}\\
\\
\ds =  \frac12\;|x|^2_{-1}+\frac12\;\sum_{k=1}^\infty \mu_k^2\;\E\int_0^t |X(s)e_k|^2_{-1}ds,
\end{array} 
\end{equation}
because by \eqref{e1.3}, 
$
|(X_\lambda-X)e_k|_{-1}\le C\lambda_k|X_\lambda-X|_{-1}
$
and so by Hypothesis \ref{h1.1}(iii)
$$
\lim_{\lambda\to 0}\sum_{k=1}^\infty\mu_k^2\E\int_0^t|X_\lambda(s)e_k|^2_{-1}ds=
\sum_{k=1}^\infty\mu_k^2\E\int_0^t|X (s)e_k|^2_{-1}ds.
$$
Next letting $\lambda$ tend to zero in \eqref{e3.1}  and using \eqref{e3.8}  we see that $\P$-a.s., for all $t\in[0,T]$,
\begin{equation}
\label{e3.18}
\langle X(t),e_j   \rangle_2=\langle  x,e_j   \rangle_2+\int_0^t\langle\eta(s),\Delta e_j   \rangle_2ds+\sum_{k=1}^\infty\mu_k\int_0^t\langle   X(s)e_k,e_j\rangle _2d\beta_k(s).
\end{equation}
Note that by continuity the  $\P$-zero set does not depend on $t\in[0\textcolor{black}{,}T]$, since
$$
\sum_{k=1}^\infty\mu_k\int_0^t\langle   X(s)e_k,e_j\rangle _2d\beta_k(s)= \int_0^t\langle    e_j,\sigma(X(s))dW(s)\rangle _2 .
$$
In order to get \eqref{e3.18} we have used the fact that by \eqref{e3.14} we have
$$
\begin{array}{l}
\ds\E\left|\int_0^t\langle X_\lambda(s)e_k,e_j\rangle_2d\beta_k(s)ds-
\int_0^t\langle X(s)e_k,e_j\rangle_2d\beta_k(s)ds\right|^2
\\
\\
\ds= \E\int_0^t\langle (X_\lambda(s)-X(s))e_k,e_j\rangle^2_2 ds\le C\lambda_j^2\lambda_k^2T
| X_\lambda-X|^2_{L^2(\Omega,C([0,T];H))}
\end{array}
$$
and therefore
$$
\lim_{\lambda\to 0}\sum_{k=1}^\infty\mu_k\int_0^t\langle   X_\lambda(s)e_k,e_j\rangle _2d\beta_kds=
\sum_{k=1}^\infty\mu_k\int_0^t\langle   X(s)e_k,e_j\rangle _2d\beta_kds.
$$
 Therefore \eqref{e3.18}  follows
and this yields, via It\^o's  formula \textcolor{black}{(applied to $\langle X(t),e_j\rangle_2^2$, $t\in [0,T]$)} and summation over $j$ that
\begin{equation}
\label{e3.19}
\begin{array}{l}
\ds\frac12\;\E|X(t)|^2_{-1}+\E\int_0^t\int_\mathcal  O\eta X d\xi\;ds\
\\
\ds=\frac12\;\E|x|^2_{-1}+\frac12\;\sum_{k=1}^\infty\mu_k^2\,\E\int_0^t|X(s)e_k|^2_{-1}ds,\;\forall\;t\in [0,T].
\end{array} 
\end{equation}
Comparing  \eqref{e3.17} and  \eqref{e3.19} we get \eqref{e3.16}.
Hence $X$ is a solution to \eqref{e1.1}  as claimed.  

\textcolor{black}{To prove uniquenss we take two solutions $X^{(1)}$ and $X^{(2)}$ with corresponding $\eta^{(1)}$ and $\eta^{(2)}$. Repeating the argument above we obtain
$$
\begin{array}{l}
\ds\frac12\;\E|X^{(1)}(t)-X^{(2)}(t)|^2_{-1}\\
\\
\ds+\E\int_0^t\int_\mathcal  O(\eta^{(1)}(s)-\eta^{(2)}(s))(X^{(1)}(s)-X^{(2)}(s))d\xi ds\\
\\
\ds=\frac12\;\sum_{k=1}^\infty\mu_k^2\,\E\int_0^t|(X^{(1)}(s)-X^{(2)}(s))e_k|^2_{-1}\;ds,\quad \forall\;t\in [0,T].
\end{array} 
$$
Since, because $\Psi$ is monotone, the second term on the left is positive, by \eqref{e1.3}, Hpothesis \ref{h1.1}(iii) this implies $X^{(1)}=X^{(2)}$ by Gronwall's lemma.}

 Finally,  if $x\ge 0$ a.e. in $\mathcal  O$ we  know by \cite[Theorem  2.2]{3} that $X_\lambda\ge 0$ $\P$-a.s. and so by \eqref{e3.14}  it follows  that $X\ge 0$,  a.e in $\Omega\times (0,T)\times \mathcal  O$ as desired. This completes the proof of   Theorem \ref{t2.2}. $\Box$\bigskip

\begin{Remark}
\em Theorem \ref{t2.2} extends to any dimension $d\ge 1$ if one modifies condition {\it (iii)}
in Hypothesis \ref{h1.1} as in \cite[Condition 4.1]{3}.
 
 \end{Remark}
 
 \begin{Remark}
\em The existence part of Theorem \ref{t2.2}  remains true for stochastic porous media equations with additive noise,  i.e.
$$
dX-\Delta\Psi(X)dt=\sqrt{Q}\;dW(t),
$$
where $\Psi$ satisfies Hypothesis \ref{h1.1} and
$$
\sqrt{Q}\;dW(t)=\sum_{k=1}^\infty\mu_ke_kd\beta_k(t)
$$
with
$$\sum_{k=1}^\infty\lambda_k^{-1}\mu_k^2<+\infty.$$
The proof is exactly the same and so, it will be omitted.
 \end{Remark}
 \begin{Proposition}
\label{p3.4}
Let  $X_\lambda, \lambda\in (0,1)$, be as above, $x\in L^4(\mathcal O)$. Assume that $\Psi$ satisfies Hypothesis \ref{h1.1} with $m=1$ and for some $\delta>0$,
\begin{equation}
\label{e3.20}
(\tilde x-\tilde y)(x-y)\ge \delta(x-y)^2,\quad \forall\;(x,\tilde x), (y,\tilde y)\in \Psi.
\end{equation}
Then $X_\lambda,X \in L^2_W(0,T;L^2(\Omega,H^1_0(\mathcal O)))$
and
\begin{equation}
\label{e3.21}
\lim_{\lambda\to 0}\E|X_\lambda-X|^2_{L^2(0,T;L^2(\mathcal O))}=0.
\end{equation}
\end{Proposition}
{\bf Proof}. A simple calculation reveals that  
$$
(\Psi_\lambda(x)-\Psi_\lambda(y))(x-y)\ge \frac\delta2\;|x-y|^2,\quad \forall\;x,y\in\R
$$
for $\lambda$ sufficiently small. Then $\tilde\Psi_\lambda$ 
defined by$\tilde\Psi_\lambda(r):=\Psi_\lambda(r)-\frac\delta2 \;r,\;r\in \R,$ is increasing and so by It\^o's formula we have
\begin{equation}
\label{e3.22}
\E|X_\lambda(t)|_2^2+\frac\delta2\;\E\int_0^t|X_\lambda(s)|_{H^1_0(\mathcal O)}^2 ds\le C.
\end{equation}
As a matter of fact, we shall apply It\^o's formula not directly to equation \eqref{e3.1} but to   equation \eqref{e3.6} (cf. the proof of Lemma \ref{l3.1} to obtain \eqref{e3.7}). Thus we get
$$
\frac12\;\E|X^\varepsilon_\lambda(t)|_2^2+\E\int_0^t\langle (A_\lambda)_\varepsilon  X^\varepsilon_\lambda(s),X^\varepsilon_\lambda(s)\rangle_2 ds\le
\frac12\;|x|_2^2+C\E\int_0^t|X^\varepsilon_\lambda(s)|_2^2 ds.
$$
Next we have
$$
\langle (A_\lambda)_\varepsilon  X^\varepsilon_\lambda ,X^\varepsilon_\lambda \rangle_2 =\langle A_\lambda(1+\varepsilon A_\lambda)^{-1}  X^\varepsilon_\lambda, (1+\varepsilon A_\lambda)^{-1}X^\varepsilon_\lambda \rangle_2+\varepsilon|(A_\lambda)_\varepsilon  X^\varepsilon_\lambda|^2_2.
$$
Taking into account that $A_\lambda=\Delta(\Psi_\lambda+\lambda I)$ and that $r\to\Psi_\lambda(r)-\delta r/2$ is monotonically increasing we get
$$
\langle (A_\lambda)_\varepsilon  X^\varepsilon_\lambda ,X^\varepsilon_\lambda \rangle_2\ge \frac\delta2\;\int_\mathcal O |\nabla(1+\varepsilon A_\lambda)^{-1}X^\varepsilon_\lambda|^2d\xi+\varepsilon |(A_\lambda)_\varepsilon  X^\varepsilon_\lambda|^2_2.
$$
Hence
$$
\E\int_0^t|(1+\varepsilon A_\lambda)^{-1}X^\varepsilon_\lambda(s)|_{H^1_0(\mathcal O )}^2 ds\le C
$$
and letting $\varepsilon\to 0$ we get \eqref{e3.22} and the first assertion (taking also into
account \eqref{e3.5}).

To prove the second part we note that
$$
\begin{array}{l}
\ds d(X_\lambda-X_\mu)-\Delta[\tilde\Psi_\lambda(X_\lambda)-
\tilde\Psi_\mu(X_\mu)+\lambda X_\lambda-\mu X_\mu+\frac12\;\delta\,(X_\lambda-  X_\mu)]dt
\\
\\
=(\sigma(X_\lambda)-\sigma(X_\mu))dW.
\end{array}
$$
Hence exactly the same arguments to derive \eqref{e3.11} lead to
$$
\begin{array}{l}
\ds \frac12\;|X_\lambda(t)-X_\mu(t)|^2_{-1} e^{-\alpha t}+\frac\delta2\;\int_0^t|X_\lambda(s)-X_\mu(s)|^2_{2} e^{-\alpha s}ds
\\
\\
\ds \le C\max\{\lambda,\mu\}\int_0^t \Big(|\Psi_\lambda(X_\lambda(s))|_2^2+|\Psi_\mu(X_\mu(s))|_2^2\\
\\
\ds\hspace{20mm}+|X_\lambda(s)|^2_2+|X_\mu(s)|^2_{2}\Big)e^{-\alpha s}  ds+M_{\lambda,\mu}(t),
\end{array}
$$
for $\alpha$ large enough and $\lambda,\mu\in (0,1)$, $t\in [0,T]$.
Since $m=1$, we have $|\Psi_\lambda(x)|\le C(1+|x|)$
for all $x\in \R$, $\lambda\in (0,1)$, hence taking expectation we get
$$
\frac\delta2\;\E\int_0^t|X_\lambda(s)-X_\mu(s)|^2_{2}ds\le C
\max\{\lambda,\mu\}\E\int_0^t(|X_\lambda(s)|^2+|X_\mu(s)|^2)ds.
$$
By Lemma \ref{l3.1} with $p=2$ and \eqref{e3.8} this implies \eqref{e3.21}. $\Box$\bigskip

Besides Hypothesis \ref{h1.1}, we shall now assume   the following
 {\it \begin{enumerate}
\item[{\it (iv) }]  $\Psi(r)=\rho\;\mbox{\rm sign}\;r+\widetilde{\Psi}(r)$,  for   $r\in \R,$
where $\rho>0$, $\widetilde{\Psi}:\R\to\R$ is  Lipschitzian, $\widetilde{\Psi}\in C^1(\R\setminus\{0\})$  and  for some $\delta>0$ it  satisfies $\widetilde{\Psi}'(r)\ge \delta $ for all $r\in \R\setminus\{0\}$.
\end{enumerate}}
Here the  signum is defined by \eqref{e1.6}.  

Below we shall use an approximation to $\Psi$ which is slightly different from $\Psi_\lambda$ defined before. Namely, below we consider
$$
\Psi_\lambda(r):=\rho\;(\mbox{\rm sign})_\lambda(r)+\widetilde{\Psi}(r),\quad r\in \R,
$$
where $(\mbox{\rm sign})_\lambda$  is the  Yosida approximation of the sign, i.e. $$
 (\mbox{\rm sign})_\lambda(r):=
\left\{\begin{array}{l}
\protect
1\quad\mbox{\rm if}\;  r>\lambda\\
\frac{r}\lambda\quad\mbox{\rm if}\;  r\in[-\lambda,\lambda]\\
-1\quad\mbox{\rm if}\;  r<-\lambda.
\end{array}\right.
$$
We shall use the   symbol $\Psi_\lambda$ also for this approximation  and denote also by $X_\lambda$  the corresponding solution  of \eqref{e3.1}.
This approximation in the special case of condition (iv) is much more convenient. We emphasize that all previous results   remain true for this modified approximation. The proofs are the same and some parts even simplify. We therefore shall use  all previous results  for $\Psi_\lambda$ and $X_\lambda$ as above without further notice.

\begin{Proposition}
\label{p3.5}
The solutions $X_\lambda$  to \eqref{e3.1} and $X$ to   \eqref{e1.1}
satisfy all conditions of Proposition \ref{p3.4} and in addition
$$
\E\int_0^T\int_\mathcal O|\nabla(\mbox{\rm sign})_\lambda(X_\lambda)|^2 d\xi dt\le C,\quad\forall\;\lambda>0
$$
and consequently $\eta\in L^2_W(0,T;L^2(\Omega;H^1_0(\mathcal O))$.
\end{Proposition}
{\bf Proof}. We set
$$
g_\lambda(r):=\int_0^r(\mbox{\rm sign})_\lambda(s)ds,\quad r\in \R,
$$
and choose $\varphi_\lambda\in C^2(\R)$ such that
\begin{enumerate}
\item[(i)] $\varphi_\lambda(0)=0$.
\item[(ii)] $\varphi'_\lambda(r)=\frac{r}\lambda$ for $|r|\le \lambda$,
$\varphi'_\lambda(r)=1+\lambda$ for $r\ge 2\lambda$,
$\varphi'_\lambda(r)=-1-\lambda$ for $r\le -2\lambda$.

\item[(iii)] $0\le \varphi''_\lambda(r) \le \frac{C}{\lambda}$ for   all $r\in \R$.
\end{enumerate}
It is easily seen that such a function exists and can be   constructed
simply by smoothing the function (sign)$_\lambda$.
Let us denote the resulting function by $f_\lambda$. Then define
$$
\varphi_\lambda(r):=\int_0^rf_\lambda(s)ds,\quad r\in \R
$$
As mentioned above the arguments of the previous proofs extends to the present situation  in order to prove that $\{X_\lambda\}$ is convergent to the solution $X$ to  \eqref{e1.1}.

Now we shall apply It\^o's formula to equation \eqref{e3.1} (or, more exactly,  to \eqref{e3.6}  and then let $\varepsilon\to 0$ as in the proof of Proposition \ref{p3.4})
with $\Psi_\lambda$ defined as above and to the function $\int_\mathcal O\varphi_\lambda(X_\lambda)d\xi$. 

Arguing as in the proof of Lemma \ref{l3.1} to obtain \eqref{e3.7}, we get (recall that $X_\lambda(t)\in H^1_0(\mathcal O)$),
$$
\begin{array}{l}
\ds\E\int_\mathcal O\varphi_\lambda(X_\lambda(t))d\xi-\E\int_0^t\langle \Delta(\mbox{\rm sign})_\lambda(X_\lambda(s))+ \Delta\widetilde{\Psi}(X_\lambda(s)), 
\varphi'_\lambda(X_\lambda(s))\rangle_2\;ds
\\
\\
\ds \le\int_\mathcal O\varphi_\lambda(x)d\xi+C\sum_{k=1}^\infty\mu_k^2 \E\int_0^t\int_\mathcal  O\varphi''_\lambda(X_\lambda(s))|X_\lambda(s)e_k|^2   d\xi ds
\\
\\
\ds \le \int_\mathcal O\varphi_\lambda(x)d\xi+4\lambda C\;\sum_{k=1}^\infty\mu_k^2 \lambda_k^2\E\int_0^t\int_\mathcal  O 1_\lambda(s,\xi)|e_k|^2   d\xi ds,
\end{array}
$$
where $1_\lambda$ is the characteristic function of the set $\{(s,\xi):\;
0\le |X_\lambda(s,\xi)|\le 2\lambda\}$.

Concerning the first  line
we note that, since $\varphi'_\lambda$ and $\tilde \Psi$ are monotonically increasing
while
as seen earlier $X_\lambda(t)\in H^1_0(\mathcal O)$, we have by the Green formula that
$$ \langle \Delta \widetilde\Psi(X_\lambda), \varphi_\lambda'(X_\lambda)\rangle_2 =-\int_\mathcal O \widetilde\Psi  '(X_\lambda)\varphi_\lambda''(X_\lambda)|\nabla X_\lambda|^2d\xi\le 0.$$
This yields
$$
\E\int_0^T\int_\mathcal O\langle \nabla(\mbox{\rm sign})_\lambda(X_\lambda),\nabla\varphi'_\lambda(X_\lambda)\rangle_2\;d\xi ds\le C,\quad \forall\;\lambda\in (0,1).
$$ 
Taking into account that
$$
-\langle \Delta(\mbox{\rm sign})_\lambda(X_\lambda),\varphi'_\lambda(X_\lambda)\rangle_2
=\langle \nabla(\mbox{\rm sign})_\lambda(X_\lambda),\nabla\varphi'_\lambda(X_\lambda)\rangle_2\ge 0,\quad\mbox{\rm a.e.}
$$
and that $\nabla\varphi'_\lambda(X_\lambda)=\frac1\lambda\;
\nabla X_\lambda$ on $\{(s,\xi):\;
|X_\lambda(s,\xi)|< \lambda\}$ we get
$$
\E\int_0^T\int_\mathcal O|\nabla(\mbox{\rm sign})_\lambda(X_\lambda)|^2
 d\xi ds\le C,\quad \forall\;\lambda\in (0,1),
$$
because $\nabla(\mbox{\rm sign})_\lambda(X_\lambda)=\frac1\lambda
\nabla(X_\lambda)$ if $|X_\lambda)|< \lambda$ and
 $\nabla(\mbox{\rm sign})_\lambda(X_\lambda)=0$ if $|X_\lambda)|\ge \lambda$.
 
 Then we get the desired estimate and since also by \eqref{e3.22}
 $$
 \E\int_0^T\int_\mathcal O|\nabla \widetilde{\Psi}(X_\lambda)|^2 d\xi ds\le C,\quad\forall\;\lambda\in (0,1)
 $$
 and $(\mbox{\rm sign})_\lambda(X_\lambda)+\widetilde{\Psi}(X_\lambda)\to \eta$ weakly  in $L^2(\Omega\times (0,T)\times\mathcal O)$ as $\lambda\to 0$ we infer that $\eta\in L^2_W(0,T;L^2(\Omega;H^1_0(\mathcal O))$ as claimed. $\Box$

 \section{Extinction in finite time and self-organized criticality}
 In this section we shall prove   a finite extinction property for solutions of  \eqref{e1.1} in $1$-$D$ for a special density dependent diffusion coefficient function $\Psi$. 
 However, Lemma \ref{l4.1} below can be proved without restriction on dimension. So,   for the moment we remain in our general framework.

 For simplicity we choose   the Wiener process
 \begin{equation}
\label{e4.1}
W(t)=\sum_{k=1}^N\mu_k e_k\beta_k(t),\quad t\ge 0,
\end{equation}
where $N\in \N$. 
 
 Besides Hypothesis \ref{h1.1}, we shall assume   Hypothesis (iv) (see  page 16), i.e. 
 {\it \begin{enumerate}
\item[{\it (iv) }]  $\Psi(r)=\rho\;\mbox{\rm sign}\;r+\widetilde{\Psi}(r)$,  for   $r\in \R,$
where $\rho>0$, $\widetilde{\Psi}:\R\to\R$ is  Lipschitzian, $\widetilde{\Psi}\in C^1(\R\setminus\{0\})$  and  for some $\delta>0$ it  satisfies $\widetilde{\Psi}'(r)\ge \delta $ for all $r\in \R\setminus\{0\}$.
\end{enumerate}}
Here the  signum is defined by \eqref{e1.6}.

Now let $\tau$ be the stopping time
$$ 
\tau=\inf\{t\ge 0:\;|X(t,x)|_{-1}=0\},
$$ 
 where $X(t,x), t\ge 0,$ is the solution to \eqref{e1.1}   given by Theorem \ref{t2.2} for $x\in L^p(\mathcal O)$, $p\ge \max\{4,2m\}.$
 \begin{Lemma}
\label{l4.1}
Under assumptions {\it (i)-(iv) } we have
$$
X(t,x)=0,\quad\mbox{\rm for} \; t\ge \tau,\;\;\P\mbox{\rm -a.s.}.
$$
\end{Lemma}
{\bf Proof}. Set $A=-\Delta$, $D(A)=H^2(\mathcal O)\cap H^1_0(\mathcal O).$
Define $\mu: [0,T]\times \Omega\to C^2_b(\mathcal O;\R)$ by
$$
\mu(t):=-\sum_{k=1}^N\mu_ke_k\beta_k(t),\quad t\in [0,T],
$$
and $\tilde\mu: [0,T]\to C^2_b(\mathcal O;\R)$ by
$$
\tilde\mu:=\sum_{k=1}^N\mu^2_ke^2_k.
$$
Define
$$
Y(t)=e^{\mu(t)}X(t),\quad t\ge 0.
$$
Let $D(A)$ be equipped with the graph norm of $A$ and let $D(A)'$ be its dual space, hence
\begin{equation}
\label{e4.2}
D(A)\subset H^1_0(\mathcal O)\subset L^2(\mathcal O)\subset H^{-1}(\mathcal O)\subset D(A)'.
\end{equation}
It is easy to see that for all $\omega\in \Omega$, $ t\in [0,T]$ the function
 $e^{\mu(t,\omega)}$ is a multiplier both in $D(A)$ and in $H$, hence
$e^{\mu(t,\omega)}\Delta z\in D(A)'$ is well defined for all $z\in L^2(\mathcal O)$ and $Y(t)\in H$.\bigskip

{\bf Claim}. We have
\begin{equation}
\label{e4.3}
Y(t)=x+\int_0^te^{\mu(s)}\Delta\eta(s) ds-\frac12\;\int_0^t\tilde\mu Y(s)ds,\quad t\in [0,T],
\end{equation}
where the fist integral on the right hand side is a Bochner  integral in 
$D(A)'$, the second by \eqref{e3.8} is one in $L^p(\mathcal O)\subset L^2(\mathcal O)$. In particular a posteriori the first integal is in $H$, continuous in $H$ as a function of $t\in [0,T]$, $\P$-a.s.\bigskip

{\bf Proof of the Claim}. Let $\varphi\in D(A)$. As before we shall use $\langle \cdot, \cdot  \rangle_2$ also for the extended dualizations with pivot space
$L^2(\mathcal O)$ as the ones in \eqref{e4.2}.Then for $t\in [0,T]$
$$
\langle \varphi, e^{\mu(t)}X(t)  \rangle_2=\sum_{j=1}^\infty
\langle e_j, e^{\mu(t)}\varphi\rangle_2\;\langle e_j, X(t)\rangle_2
$$
Furthermore, we have by It\^o's formula for all $\xi\in \mathcal O$
$$
e^{\mu(t,\xi)}=1+\int_0^te^{\mu(s,\xi)}d\mu(s,\xi)+\frac12\;\int_0^te^{\mu(s,\xi)}\tilde\mu(\xi)ds.
$$
Now fix $j\in \N$. Then by the stochastic Fubini Theorem
$$
\begin{array}{lll}
\ds \langle e_j, e^{\mu(t)}\varphi\rangle_2&=&\ds\langle e_j,  \varphi\rangle_2
-\sum_{k=1}^N\mu_k\int_0^t\langle e_j, e_ke^{\mu(s)}\varphi\rangle_2d\beta_k(s)
\\
\\
&&\ds+\frac12\;\int_0^t\langle e_j,  \tilde\mu e^{\mu(s)}\varphi\rangle_2ds,\quad t\in [0,T].
\end{array}
$$
By It\^o's product rule and \eqref{e3.18} we hence obtain 
$$
\begin{array}{l}
\ds \langle e_j, e^{\mu(t)}\varphi\rangle_2\;\langle e_j, X(t)\rangle_2
=\langle e_j,  \varphi\rangle_2\;\langle e_j, x\rangle_2
\\
\\
\ds+\int_0^t\langle e_j,    e^{\mu(s)}\varphi\rangle_2\;\langle \Delta e_j, \eta(s)\rangle_2\;ds\\
\\
\ds+\sum_{k=1}^N\mu_k\int_0^t\langle e_j,   e^{\mu(s)}\varphi\rangle_2\;\langle   e_j, X(s)e_k\rangle_2\;d\beta_k(s)\\
\\
\ds+\frac12\;\int_0^t\langle e_j,  X(s)  \rangle_2\;\langle  e_j,  \tilde\mu \;e^{\mu(s)}\varphi\rangle_2\;ds\\
\\
\ds-\sum_{k=1}^N\mu_k\int_0^t\langle e_j,   X(s)\rangle_2\;\langle   e_j, e_ke^{\mu(s)}\varphi\rangle_2\;d\beta_k(s)
\\
\\
\ds-\sum_{k=1}^N\mu^2_k\int_0^t\langle e_j,   e_ke^{\mu(s)}\varphi\rangle_2\;\langle   e_j, X(s)e_k\rangle_2\;d\beta_k(s).
\end{array}
$$
After summing over $j\in \N$ the two stochastic terms cancel and the claim follows since $\varphi\in D(A)$ was arbitrary.

Below we work for $\P$-a.s. $\omega\in \Omega$, $\omega$ fixed. Hence all constants $C$ appearing below may depend on $\omega$.

Consider the solution $X_\lambda\in L^2_W(0,T;L^2(\Omega,H^1_0(\mathcal O)))$ to equation \eqref{e3.1}. By Proposition \ref{p3.4} we have
$$
\lim_{\lambda\to 0}\E|X_\lambda-X|^2_{L^2(0,T;L^2(\mathcal O))}=0
$$
and $\Psi_\lambda(X_\lambda)\in L^2_W(0,T;L^2(\Omega,H^1_0(\mathcal O)))$ because $\Psi_\lambda$ is Lipschitz.

On the other hand,we have as in \eqref{e4.3} for $Y_\lambda=e^\mu X_\lambda$
\begin{equation}
\label{e4.4}
\frac{dY_\lambda(t)}{dt}=e^{\mu(t)}\Delta \eta_\lambda(t)-\frac12\;\tilde\mu(t)Y_\lambda(t),\quad \forall\;t\ge 0
\end{equation}
where
$$
\eta_\lambda(t)=\Psi_\lambda(X_\lambda(t)) \in H^1_0(\mathcal O).
$$
It follows by  \eqref{e3.21} that
\begin{equation}
\label{e4.5}
\lim_{\lambda\to 0}\E|Y_\lambda-Y|^2_{L^2(0,T;L^2(\mathcal O))}=0
\end{equation}
and therefore for some sequence $\lambda_n\to 0$
\begin{equation}
\label{e4.6}
\lim_{n\to \infty}|Y_{\lambda_n}-Y|_{L^2(0,T;L^2(\mathcal O))}=0\quad\mbox{\rm a.e. on}\;\Omega. 
\end{equation}
Below we simple write $\lambda$ instead of $\lambda_n$. Next we have by \eqref{e4.4}  that
\begin{equation}
\label{e4.7}
\left<\frac{dY_\lambda(t)}{dt},Y_\lambda(t)\right>_2=
\left<\eta_\lambda(t),\Delta(e^{\mu(t)}Y_\lambda(t))\right>_2-\frac12\;\langle\tilde\mu(t)Y_\lambda(t),Y_\lambda(t)  \rangle_2\quad\mbox{\rm a.e.}\;t\in [0,T].
\end{equation}

Also we have (for simplicity we take $\rho=1$)
$$
\begin{array}{l}
\ds \langle\eta_\lambda(t),\Delta(e^{\mu(t)}Y_\lambda(t))   \rangle_2
\\
\\
\ds=\langle(\mbox{\rm sign})_\lambda\;(e^{-\mu(t)}Y_\lambda(t))+\widetilde{\Psi}(e^{-\mu(t)}Y_\lambda(t)) ,\Delta(e^{\mu(t)}Y_\lambda(t))   \rangle_2\\
\\
\ds
=-\int_\mathcal O(\nabla (\mbox{\rm sign})_\lambda\;(e^{-\mu(t)}Y_\lambda(t)),\nabla(e^{\mu(t)}Y_\lambda(t))) d\xi\\
\\
\ds -\int_\mathcal O
\widetilde{\Psi}'(e^{-\mu(t)}Y_\lambda(t))
(\nabla  (e^{-\mu(t)}Y_\lambda(t)),\nabla(e^{\mu(t)}Y_\lambda(t))) d\xi\\
\\
\ds = -\frac1\lambda\;\int_\mathcal O(|\nabla Y_\lambda(t)|^2-|Y_\lambda(t)|^2 
\;|\nabla\mu(t)|^2)1_\lambda(t,\xi) d\xi\\
\\
\ds   - \int_\mathcal O\widetilde{\Psi}'(e^{-\mu(t)}Y_\lambda(t))(|\nabla Y_\lambda(t)|^2-|Y_\lambda(t)|^2 
\;|\nabla\mu(t)|^2)  d\xi,
\end{array}
$$
because for $y\in H^1_0(\mathcal O)$ 
$$
\nabla\;(\mbox{\rm sign})_\lambda\;(y)=
\left\{\begin{array}{l}
\protect
0,\quad\mbox{\rm on}\;\{y\notin (-\lambda,\lambda)\},\\\\
\frac1\lambda\;\nabla y,\quad\mbox{\rm on}\;\{y\in (-\lambda,\lambda)\}.
\end{array}\right.
$$
(Here $1_\lambda$ is the characteristic function of $\{(\xi,t)\in \mathcal O\times [0,T]:\;|e^{-\mu(t,\xi)}Y_\lambda(t,\xi))|<\lambda   \}$ and $(\cdot,\cdot)$ is the euclidean scalar product in $\R^n$.)
Since $\widetilde{\Psi}'\ge \delta$ and $\widetilde{\Psi}'\in L^\infty(\R)$, $\mu\in C([0,T]\times\mathcal O)$ this yields
\begin{equation}
\label{e4.7'}
\langle\eta_\lambda(t),\Delta(e^{\mu(t)}Y_\lambda(t))   \rangle_2
\le C\left(|Y_\lambda(t)|^2_2+\lambda\right).
\end{equation}
Hence \eqref{e4.7} and Gronwall's lemma imply
$$
|Y_\lambda(t)|^2_2\le e^{C(t-s)}\left(|Y_\lambda(s)|^2_2+C\lambda T\right)\quad\mbox{\rm a.e.}\;t>s.
$$
Now taking into account \eqref{e4.6} and letting $\lambda\to 0$  we get
\begin{equation}
\label{e4.8}
|Y(t)|^2_2\le e^{C(t-s)}|Y(s)|^2_2\quad\mbox{\rm a.e.}\;t>s.
\end{equation}
If  $Y(\cdot)$ is  
$L^2(\mathcal O)$-continuous   then  \eqref{e4.8} holds for all $s,t\in [0,T]$, $t\ge s$. Taking in \eqref{e4.8} $s=\tau\wedge T$ we get $Y(t)=0$ for all $t\ge \tau\wedge T$ and since $T>0$ was arbitrary for all $t\ge \tau$ as claimed. So,
we have to prove that $Y$ is $L^2(\mathcal O)$-continuous on $[0,T]$. For this we recall that by Proposition \ref{p3.5} we have
 \begin{equation}
\label{e4.9}
e^{\mu}\;\eta\in L^2(0,T;H^1_0(\mathcal O)),\quad \P\mbox{\rm -a.s.}.
\end{equation} 
 Then by equation \eqref{e4.3} we have $\frac{dY}{dt} \in L^2(0,T;H^{-1}(\mathcal O))$ and so, since $Y\in L^2(0,T;H^1_0(\mathcal O))$ $\P$-a.s. by
 Proposition \ref{p3.4}, by a well known interpolation result (see e.g. \cite{0}), we conclude that
 $Y  \in C([0,T];L^2(\mathcal O))$.  
 This concludes the proof of Lemma \ref{l4.1}. $\Box$

  \bigskip

For proving our extinction result we need $\mathcal O\subset \R$, i.e. $d=1$. To be more specific let $\mathcal O=(0,\pi)$. Then
$e_k(\xi)=\sqrt{\frac2\pi}\;\sin k\xi,\quad \xi\in [0,\pi]$, $\lambda_k=k^2$ and $L^1(0,\pi)\subset H$ continuously, so
\begin{equation}
\label{e4.11}
\gamma=\inf\left\{ \frac{|x|_{L^1}}{|x|_{-1}}:\;x\in L^1(0,\pi)  \right\}>0.
\end{equation}

\begin{Theorem}
\label{t4.2}
Let  $x\in L^p(0,\pi),$  $p\ge \max\{2m,4\},$ be such that
$$
|x|_{-1}< C_N^{-1} \rho\gamma,
$$
where
\begin{equation}
\label{e4.12}
C_N:=\frac\pi{4}\;\sum_{k=1}^N(1+k)^2\mu_k^2.
\end{equation}
Then, for each $n\in \N$,
\begin{equation}
\label{e4.13}
\P(\tau\le n)\ge 1-\textcolor{black}{\frac{|x|_{-1}}{\rho\gamma}\;\left(\int_0^ne^{-C_Ns}ds\right)^{-1}},
\end{equation}
where by Lemma \ref{l4.1} we have
$$
\tau(\omega)=\sup\{t\ge 0:\; |X(t,x)|_{-1}>0\}.
$$
\end{Theorem}

\noindent{\bf Proof}.  By condition {\it (iv)} we see that
\begin{equation}
\label{e4.15}
r\Psi(r)\ge \rho|r|,\quad \forall\;r\in \R.
\end{equation}
Consider the solution $X_\lambda\in L^2_W(0,T;L^2(\Omega;H^1_0(0,\pi)))$ to equation \eqref{e3.1}. Then by first applying Krylov-Rozovskii's It\^o formula (cf.\cite[Theorem I.3.1]{KR} or e.g.
\cite[Theorem 4.2.5]{PR}) and then the classical It\^o formula  to the real valued semi-martingale $|X_\lambda(t)|^2_{-1}, t\in [0,T],$
and the function
$$
\varphi_\varepsilon(r)=(r+\varepsilon^2)^{1/2},\quad r\in \R,
$$
we find
\begin{equation}
\label{e4.16}
\begin{array}{l}
d\varphi_\varepsilon(|X_\lambda(t)|^2_{-1})+(|X_\lambda(t)|^2_{-1}+\varepsilon^2)^{-1/2}\langle X_\lambda(t),\Psi_\lambda(X_\lambda(t))    \rangle_2 dt \\
\\
\ds=\frac12\;\sum_{k=1}^N\mu_k^2
\frac{|X_\lambda(t)e_k|^2_{-1}(|X_\lambda(t)|^2_{-1}+\varepsilon^2)-|\langle X_\lambda(t)e_k,X_\lambda(t)   \rangle_{-1}|^2)}{(|X_\lambda(t)|_{-1}^2+\varepsilon^2)^{3/2}}\;dt\\
\\
+\langle\sigma(X_\lambda(t))dW(t),\varphi'_\varepsilon(|X_\lambda(t)|^2_{-1})X_\lambda(t)   \rangle
\\\\
\ds\le\frac12\;\sum_{k=1}^N\mu_k^2
\frac{|X_\lambda(t)e_k|^2_{-1}}{(|X_\lambda(t)|_{-1}^2+\varepsilon^2)^{1/2}}dt+\langle\sigma(X_\lambda(t))dW(t),\varphi'_\varepsilon(|X_\lambda(t)|^2_{-1})X_\lambda(t) \rangle
\\
\\
\ds\le C_N
\frac{|X_\lambda(t)|^2_{-1}}{(|X_\lambda(t)|_{-1}^2+\varepsilon^2)^{1/2}}\;dt+\textcolor{black}{2}\langle\sigma(X_\lambda(t))dW(t),\varphi'_\varepsilon(|X_\lambda(t)|^2_{-1})X_\lambda(t)  \rangle.
\end{array} 
\end{equation}
Here $C_N$ is given by \eqref{e4.12}   and
$$
\sigma(X_\lambda(t))dW(t)=\sum_{k=1}^N\mu_kX_\lambda(t)e_kd\beta_k(t).
$$
Integrating over $t$ and letting $\lambda\to 0$ we see that the right hand side of \eqref{e4.16} converges to  the right hand side of \eqref{e4.17} below in $L^2(\Omega;C([0,T];H))$. But by \eqref{e3.5},\eqref{e3.8}, \eqref{e3.12},
\eqref{e3.13} and by Proposition \ref{p3.4} the same is true for the left hand side with limit
$$
\varphi_\varepsilon(|X(t)|^2_{-1})-\varphi_\varepsilon(|x|^2_{-1})
+\int_0^t\int_\mathcal O \frac{X(s)}{(|X(s)|^2_{-1}+\varepsilon)^{1/2}}\;\eta(s)d\xi ds.
$$
Taking into account    \eqref{e2.2} and \eqref{e4.15} we altogether obtain
$$
\begin{array}{l}
\ds d\varphi_\varepsilon(|X(t)|^2_{-1})+\rho\frac{|X(t)|_{L^1(0,\pi)}}{(|X(t)|_{-1}^2+\varepsilon^2)^{1/2}}dt\\
\\
\ds\le C_N
\frac{|X(t)|^2_{-1}}{(|X(t)|_{-1}^2+\varepsilon^2)^{1/2}}\;dt
+\textcolor{black}{2}\langle\sigma(X(t))dW(t),\varphi'_\varepsilon(|X(t)|^2_{-1})X(t)   \rangle.
\end{array} 
$$
Consequently by Lemma \ref{l4.1} for all $t\ge 0$
\begin{equation}
\label{e4.17}
\begin{array}{l}
\ds\varphi_\varepsilon(\textcolor{black}{|X(t)}|^2_{-1})+\gamma \rho\int_0^{t\wedge\tau}\frac{|X(s)|_{-1}}{(|X(s)|_{-1}^2+\varepsilon^2)^{1/2}}ds\\
\\
\ds\le\varphi_\varepsilon(|x|^2_{-1})+C_N
\int_0^{t\wedge\tau}  \frac{|X(s)|^2_{-1}}{(|X(s)|_{-1}^2+\varepsilon^2)^{1/2}}ds\\
\\
\ds+\textcolor{black}{2}\int_0^{t\wedge\tau}\langle\sigma(X(s))dW(s),\varphi'_\varepsilon(|X(s)|^2_{-1})X(s) \rangle,
\quad\P\mbox{\rm -a.s.},
\end{array} 
\end{equation}
where $\gamma$ is defined by  \eqref{e4.4}.

Clearly, we have  
$$
\lim_{\varepsilon\to 0}\int_0^{t\wedge\tau}\frac{|X(s)|_{-1}}{(|X(s)|^2_{-1}+\varepsilon^2)^{1/2}}\;ds=t\wedge\tau,\quad   \quad\P\mbox{\rm -a.s.}.
$$
Now, letting $\varepsilon$ tend to zero we get
\begin{equation}
\label{e4.18}
\begin{array}{l}
\ds |X(t)|_{-1}+\gamma\rho (t\wedge\tau)\le|x|_{-1}+C_N
\int_0^t  |X(s)|_{-1}ds\\
\\
\ds+\int_0^t1_{[0,\tau]}(s)\langle \sigma(X(s))dW(s),X(s)|X(s)|^{-1}_{-1}  \rangle
\quad\P\mbox{\rm -a.s.}
\end{array} 
\end{equation}
Hence by a standard comparison result
$$
\begin{array}{l}
\ds|X(t)|_{-1}+\rho\gamma\int_0^te^{C_N(t-s)}1_{[0,\tau]}(s)ds\le e^{C_Nt}|x|_{-1}
\\
\\
\ds+
\int_0^te^{C_N(t-s)}1_{[0,\tau]}(s)\langle \sigma(X(s))dW(s),X(s)|X(s)|^{-1}_{-1}  \rangle.
\end{array}
$$
Taking expectation and multiplying by $(\rho\gamma)^{-1}e^{-C_N t}$, we obtain
$$
\int_0^te^{-C_Ns}\P(\tau>s)ds\le \frac{|x|_{-1}}{\rho\gamma}.
$$
Writing $\P(\tau>s)=1-\P(\tau\le s)$ we deduce that
$$
\P(\tau\le t)\ge 1-\textcolor{black}{\frac{|x|_{-1}}{\rho\gamma}\;\left(\int_0^te^{-C_Ns}ds\right)^{-1}}
$$
and \eqref{e4.13} follows. $\Box$

In particular Theorem \ref{t4.2} applies to self-organized criticality stochastic models \eqref{e1.7}
\begin{equation}
\label{e4.19}
\left\{\begin{array}{l}
dX(t)-\Delta(\rho\;\mbox{\rm sign}\;(X(t)-x_c)+\widetilde{\Psi}(X(t)-x_c))dt\\
\\
\hspace{20mm}\ds\ni \sigma (X(t)-x_c)\sum_{k=1}^N\mu_ke_kd\beta_k,\quad t\ge 0, \\
\\
\rho\;\mbox{\rm sign}\;(X(t)-x_c)+\widetilde{\Psi}(X(t)-x_c)\ni 0,\quad\mbox{\rm on}\;\partial [0,\pi],\\
\\
X(0,x)=x.
\end{array}\right.
\end{equation}
Here the function $\widetilde{\Psi}$ is as in assumption {\it (iv)} and $x_c\in \R$.
\begin{Corollary}
\label{c4.3}
Assume that
$$
|x-x_c|_{-1}<\rho\gamma C_N^{-1},
$$
where $C_N$ is as in \eqref{e4.12} and $\gamma$ as in \eqref{e4.11}.
Then for each $n\in \N$
\begin{equation}
\label{e4.20}
\P(\tau_c\le n)\ge 1-\textcolor{black}{\frac{|x-x_c|_{-1}}{\rho\gamma}\;\left(\int_0^ne^{-C_Ns}ds\right)^{-1}}
,
\end{equation}
where
$$
\tau_c=\inf\{t\ge 0:\;|X(t)-x_c|_{-1} =0  \}=\sup\{t\ge 0:\;
|X(t)-x_c|_{-1}>0\}.
$$
and $X=X(t,x)$  is the solution to \eqref{e4.19} in the sense of Definition 2.1.\end{Corollary}

We note that equation \eqref{e1.7} reduces to  \eqref{e4.19}  by shifting the Heavside function with $x_c$.

One must notice that if $x>x_c,$ i.e.   if the initial state is in the supercritical region  then by positivity result in Theorem \ref{t2.2} we have
$X(t)\ge x_c$, $\P$-a.s. for all $t\ge 0$. This means that the state remains in 
the supercritical-critical region for all  the time. However, by \eqref{e4.20} if $\frac{C_N|x|_{-1}}{\rho\gamma}$ is small, it reaches the critical state $x_c$ with high probability  in a finite time i.e. the supercritical-critical region is completely absorbed by the critical one in a finite time.
\begin{Remark}
\label{r4.4}
\em 
\textcolor{black}{Let us consider the deterministic case. Then \eqref{e4.13} implies that $\tau\le n$ if $n>\frac{|x|_{-1}}{\rho\gamma};$ so $\tau\le \frac{|x|_{-1}}{\rho\gamma}$. But this, of course, also follows directly from \eqref{e4.18}, since we assume $C_N=0$.}

\end{Remark}

{\bf A\textcolor{black}{c}knowledgement.} 
We would like to thank Philippe Blanchard for introducing us to these models of self-organized criticality and the relevant literature.
This work has been supported in part by
 the CEEX Project 05 of
Romanian Minister of Research,
the DFG -International Graduate School ``Stochastics and Real World Models'',the SFB-701 and the
\textcolor{black}{BiBoS}-Research Center.', 
  the research programme ``Equazioni di
Kolmogorov'' from the Italian
``Ministero della Ricerca Scientifica e Tecnologica''
and "FCT, POCTI-219, FEDER".
 Most of the work 
was done during very pleasant visits of the first 
and second  author to the 
University of Bielefeld and the first and the third author to the SNS in Pisa.

\end{document}